\documentclass[12pt,a4paper]{article} % a4
\usepackage{amsmath}
\usepackage{amssymb}
\usepackage{times} 
\newcommand{\nc}{\newcommand} 
\nc{\bb}{\bigskip}
\nc{\cl}{\centerline}
\nc{\e}{\varepsilon} 
\nc{\ind}{\hskip 2em\relax}
\nc{\la}{\langle}\nc{\ra}{\rangle}
\nc{\mb}{\mathbb}
\nc{\mc}{\mathcal}
\nc{\ssi}{\Leftrightarrow}
\nc{\sx}{\mbox{x}}
\nc{\pt}{\bullet}
\nc{\then}{\Rightarrow}
\nc{\vide}{\varnothing}
\newtheorem{ttt}{Theorem}%[section]
\newtheorem{metat}[ttt]{Metatheorem}

\nc{\citeb}{\cite{bkiens}}

\begin{document}\large

\cl{\Huge There are only countably many sets}\vskip4em
%\title{There are only countably many sets}\maketitle

{\footnotesize  Comme le Dieu des philosophes, l'op\'eration de Hilbert est 
incompr\'ehensible et ne se voit pas~; mais elle gouverne tout, et ses 
manifestations sensibles \'eclatent partout.
\hfill Godement}\vskip4em

We freely use the contents of Bourbaki's book \cite{bkiens}. In particular 
$(T|x)R$ means ``$T$ replaces $x$ in $R$'' whenever $T$ is a term, $x$ a 
letter and $R$ a relation. \bb 

\begin{metat} There are only countably many sets. 
Different sets can be equal; more precisely, given any set $X$, there 
is a list $X_1,X_2$, \dots of pairwise distinct sets all equal to 
$X$; in particular any nonemp\-ty set has infinitely many elements. 
\end{metat}

{\bf Metaproof.} According to \cite{bkiens} (II.1) the words 
``set'' and ``term'' are strictly synonymous. There are only countably 
many terms. (Some of these sets are called ``uncountable''. Call them 
``funny'' if you wish. Don't be afraid of Virginia Woolf!) The terms 
$X,X\cup\vide,X\cup\vide\cup\vide,\dots$ are all different but equal. 
$\square$ 

\begin{metat} The relation $\tau_x(x=x)=\tau_x(x\not=x)$ is undecidable in 
 Bourbaki's theory $\mc{T}$.  
\end{metat}

{\bf Metaproof.} Denote $\tau_x(x\not=x)$ by $c$. For any letter $x$ and 
any relation $R$, define the terms $\sigma_x(R)$ and $\rho_x(R)$ as 
follows: if $R$ is a theorem of $\mc{T}$, then $\sigma_x(R)$ is $c$, and 
$\rho_x(R)$ is $\tau_x(x\not=c)$; otherwise $\sigma_x(R)$ and $\rho_x(R)$ 
are $\tau_x(R)$. Let $\mc{T}_\sigma$ (resp. $\mc{T}_\rho$) be the theory 
obtained from $\mc{T}$ by replacing $\tau$ by $\sigma$ (resp. $\rho$). 
Since the relation in the Metatheorem is true in $\mc{T}_\sigma$ but 
false in $\mc{T}_\rho$, it suffices to check that all the axioms of 
$\mc{T}$ hold in $\mc{T}_\sigma$ and $\mc{T}_\rho$. Let $\e$ be $\sigma$ 
or $\rho$. \bb

\nc{\cnine}{C9 (\citeb I.3.2)}
\nc{\cqt}{C14 (\citeb I.3.3)}
\nc{\cqz}{C15 (\citeb I.3.3)}
\nc{\cvq}{C24 (\citeb I.3.5)}
\nc{\ctht}{C30 (\citeb I.4.3)}
\nc{\sfive}{S5 (\citeb I.4.2)}
\nc{\ssept}{S7 (\citeb I.5.1)}

\ind {\it Verification of \ssept.} Let $R$ and $S$ be relations and $x$ a 
letter. We must show 

\begin{equation}\label{S7}
 ((\forall x)(R\ssi S))\then(\e_x(R)=\e_x(S)).
\end{equation}

We can assume that $R$ is a theorem and that $S$ is not. We obtain 
successively 
\begin{itemize}
 \item $R\nLeftrightarrow S$ by \cqz, 
 \item $(\exists x)(R\nLeftrightarrow S)$ by \sfive,
 \item $(\e_x(R)\not=\e_x(S))\then((\exists x)(R\nLeftrightarrow S))$ 
  by \cnine, 
 \item (\ref{S7}) by \cvq. 
\end{itemize}

\ind {\it Verification of the other axioms.} 
In the other axioms the Hilbert operation does not appear explicitly, but 
is hidden behind the quantifier $\exists$, which is defined by letting 
$(\exists x)R$ be $(\tau_x(R)|x)R$ if $x$ is a letter and $R$ a relation. 
Define $(\exists_\e x)R$ as $(\e_x(R)|x)R$. It is enough to show 
$(\exists_\e x)R\ssi(\exists x)R$, that is
$(\e_x(R)|x)R\ssi(\tau_x(R)|x)R$. We can assume again that 
$R$ is a theorem. In this case \sfive\ implies that $(\e_x(R)|x)R$ and 
$(\tau_x(R)|x)R$ are theorems too, and \cnine\ yields the conclusion. 
$\square$ \bb

\ind Say that a {\bf machine} is an assignment of a finite or infinite 
sequence $y=(y_1,y_2,\dots)$ of zeros and ones to each finite sequence 
$x=(x_1,x_2,\dots,x_n)$ of zeros and ones. 
Let $\la 0,1\ra$ be the set of all words in 0 and 1, let 
$\{0,1\}^{\mb{N}}$ be the set of all maps from $\mb{N}$ to 
$\{0,1\}$, where $\mb{N}$ denotes the set of positive integers, and let 
$\mc F$ be the set of all maps from $\la 0,1\ra$ into the 
union of $\la 0,1\ra$ and $\{0,1\}^{\mb{N}}$~:

$$\mc{F}=\left\{f:\la 0,1\ra\to
\la 0,1\ra\cup\{0,1\}^{\mb{N}}\right\}.$$

\bb\ind A machine $M$ shall be called a {\bf Bourbaki machine} if it is 
given by a term in the sense of \cite{bkiens}, 
i.e. if there is a term $T\in\mc{F}$ such that $T(x)_n$ and $M(x)_n$ 
are either both defined and equal or both undefined, whenever $x$ is a 
finite sequence of zeros and ones, and $n$ a positive integer. 

\begin{metat} Any Turing machine is a Bourbaki machine and vice versa. 
\end{metat}

{\bf Metaproof.} See Kleene \cite{kleene}, \S\ 68-69. $\square$ \bb 

\ind These considerations suggest the following definition 
of {\it theory} or {\it formal system} (we use the two expressions 
synonymously). For any set $X$ let $\la X\ra$ be the free monoid 
generated by $X$. Let $f:\la\mb{N}\ra\to\mb{N}$ be the bijection 
sending $(n_1,\dots,n_k)$ to $p_1^{n_1}\dots p_k^{n_k}$ where $p_i$ is 
the $i$-th prime, and $\la f\ra:\la\la\mb{N}\ra\ra\to\la\mb{N}\ra$ its 
functorial extension. Let $R\subset\la\mb{N}\ra$ and 
$P\subset\la R\ra$ be subsets. We think of $R$ as being the set 
of {\it relations} (or {\it formulas}, or {\it statements}, or 
{\it assertions}, as you like best) of the theory, and of $P$ as being the 
set of {\it proofs}. Say that $(R,P)$ is a {\bf theory} if the characteristic 
functions of $f(R)\subset\mb{N}$ and $f(\la f\ra(P))\subset\mb{N}$ are 
Bourbaki machines. \bb

\cl{\Large\bf Final remarks}\bb

%The first proof of the incompleteness of mathematics is supposed to be 
%G\"odel's in \cite{godel}. 
I would be most grateful to anybody who 
could give me an answer to the following questions. \bb

$\pt$ To which extent does the proof of G\"odel's First Incompleteness 
Theorem rely on the considerations below (taken from \cite{kleene}, \S\ 
41, p. 195)? 

\begin{quote}
{\normalsize%\small%\footnotesize
The terms: $0, 0', 0''$, \dots, which represent the particular natural 
numbers under the interpretation of the system, we call {\it numerals}, 
and we abbreviate them by the same symbols ``0'',  ``1'',  ``2'', \dots , 
respectively, as we use for natural numbers intuitively (\dots).  
Moreover, whenever we have introduced an italic letter, such as ``$x$'', 
to designate an intuitive natural number, then the corresponding bold 
italic letter ``$\pmb{x}$'' shall designate the corresponding 
numeral $0^{(x)}$, i.e. $0^{'\cdots'}$ with $x$ accents ($x\ge0$) \dots
\smallskip

Let $P(x_1,\dots,x_n)$ be an intuitive number-theoretic predicate. We say 
that $P(x_1,\dots,x_n)$ is {\it numeralwise expressible} in the formal 
system, if there is a formula $\mbox{P}(\sx_1,\dots,\sx_n)$ with no free 
variables other than the distinct variables $\sx_1,\dots,\sx_n$
such that, for each particular $n$-tuple of natural numbers 
$x_1,\dots,x_n$, \smallskip

(i) \ if $P(x_1,\dots,x_n)$ is true, then 
$\mbox{P}(\pmb{x}_1,\dots,\pmb{x}_n)$ is provable, and\smallskip

(ii) if $P(x_1,\dots,x_n)$ is false, then ``not 
$\mbox{P}(\pmb{x}_1,\dots,\pmb{x}_n)$'' is provable.}
\end{quote}

$\pt$ Who would expect any reasonable theory of mathematics to be 
complete?\bb

$\pt$ (Recall that a theory is {\it complete} if any relation without 
letters is true or false.) Is there a natural reason for making such a 
discrimination between relations without letters and relations with at 
least one letter? \bb

\ind Thank you to Bruno Blind, Fr\'ed\'eric Campana, Jean-Marie Didry 
and Fran\c ois Marque for their interest.

\footnotesize\hfill April 22, 2004\bb

\cl{Pierre-Yves Gaillard, D\'epartement de Math\'ematiques, Universit\'e 
Nancy~1, France}
%
%\cl{\href{http://www.iecn.u-nancy.fr/~gaillard/}%
%{http://www.iecn.u-nancy.fr/\~{}gaillard/}}

\end{document}